\documentclass[12pt,twoside,centertags,reqno]{article}
\usepackage{amsfonts,amssymb,amsmath,amsthm,latexsym,amsxtra,amscd}
\usepackage[english]{babel}
\usepackage{verbatim,alltt}
\usepackage{epsfig,epic,eepic}
\usepackage{fancyheadings}
\usepackage[ps,dvips,all]{xy} \SelectTips{cm}{12}
\newtheorem{lemma}{Lemma}[section]
\newtheorem{proposition}[lemma]{Proposition}
\newtheorem{theorem}[lemma]{Theorem}

\theoremstyle{definition}

\newtheorem{remark}[lemma]{Remark}

\date{}
\newcommand{\addresses}[3]{}
\newcommand{\emails}[3]{}
\newcommand{\classification}[1]{%
\renewcommand{\thefootnote}{}%
\footnotetext{\mbox{\hspace*{-12pt}
1991 {\it Mathematics Subject Classification.}}
{#1}%
}}
\newcommand{\keywords}[1]{\renewcommand{\thefootnote}{}%
\footnotetext{\mbox{\hspace*{-12pt}
{\it Key words.}} 
{#1}%
}}
\title{Exterior algebra methods for the construction  \\
of rational surfaces in the projective fourspace }
\author{Hirotachi Abo \and Frank-Olaf Schreyer}
\addresses{
Department of Mathematics, Colorado State University, 
Fort Collins, CO, 80523, USA
}
{
Saarbr\"ucken, Saarland, 66123, Germany
}
{
}
\emails{
abo@math.colostate.edu
}
{
schreyer@math.uni-sb.de
}
{
}
\newcommand{\C}{{\mathbb C}}

\newcommand{\F}{{\mathbb F}}
\newcommand{\GG}{{\mathbb G}}
\newcommand{\M}{{\mathbb M}}

\renewcommand{\P}{{\mathbb P}}
\newcommand{\Q}{{\mathbb Q}}

\newcommand{\Z}{{\mathbb Z}}

\newcommand{\LL}{{\mathbb L}}

\newcommand{\PP}{{\mathbb P}}

\newcommand{\h}{{\rm{h}}}
\newcommand{\HH}{{\rm{H}}}

\newcommand{\cO}{{\mathcal O}}

\newcommand{\s}{\mathcal}

\newcommand{\sE}{{\s E}}
\newcommand{\sF}{{\s F}}

\newcommand{\sI}{{\s I}}

\newcommand{\sX}{{\s X}}

\newcommand{\lra}{\longrightarrow}

\newcommand{\tensor}{\otimes}

\newcommand{\punkt}{\HHspace{-.3ex}\raise.15ex\HHbox to1ex{\HHuge.}}

\newlength{\br}
\newlength{\ho}
\DeclareMathOperator{\GL}{GL}

\DeclareMathOperator{\Spec}{Spec}
\DeclareMathOperator{\HHom}{Hom}

\DeclareMathOperator{\im}{im}

\DeclareMathOperator{\Coker}{Coker}

\DeclareMathOperator{\PGL}{PGL}

\DeclareMathOperator{\codim}{codim}

\newcommand{\integer}{\Z}

\newcommand{\ZZ}{\Z}

\begin{document}
\maketitle
\thispagestyle{empty}
\classification{%
%    INSERT here MSC Classification
Primary 14J10, 14J26; Secondary 14Q10
}
\renewcommand{\thefootnote}{\arabic{footnote}}%

\begin{abstract}
\noindent
%
%    INSERT here ABSTRACT of your article
The aim of this paper is to present a construction 
of smooth rational surfaces in projective fourspace
with degree $12$ and sectional genus $13$. 
The construction is based on exterior algebra methods, 
finite field searches and standard deformation theory. 
\end{abstract}

%
%    INSERT here BODY of your article
%

\section{Introduction}

This paper is dedicated to Gert-Martin Greuel on the occasion of his sixtieth
birthday. 
The use of computer algebra systems is essential for the proof of the
main result of this paper. It will become clear that without 
of computer algebra systems like  Singular and Plural
developed in Kaiserslautern we could not obtain the main 
result of this paper at all. 
We thank the group in Kaiserslautern for their excellent program.

Hartshone conjectured that only finitely many components of the Hilbert 
scheme of surfaces in $\P^4$ correspond to smooth rational surfaces. 
In 1989, this conjecture was positively solved by Ellingsrud and 
Peskine~\cite{ep}. 
The exact bound for the degree is, however, still open. This motivates
our search for 
smooth rational surfaces  
in $\P^4$. Examples of smooth  rational surfaces in $\P^4$ prior to this paper
were  known up to degree $11$, see \cite{ds}.
Our main result is the proof of existence of the following example.

 \begin{theorem}\label{main}
There exists a family of  smooth rational surfaces in $\P^4$ over $\C$ 
with $d=12$, $\pi=13$ and hyperplane class
\[
H \equiv 12L-\sum_{i_1=1}^2 4E_i-\sum_{i_2=3}^{11} 3E_i 
-\sum_{i_3=12}^{14} 2E_i-\sum_{i_4=15}^{21} E_i. 
\]
in terms of a plane model.
\end{theorem}
Abstractly these surfaces arise as the blow up of $\P^2$ in $21$ points.
$L$ and $E_i$ in the Theorem denote the class of a general line 
and the exceptional 
divisors.  

The $21$ points lie in special position due to the
fact that we need $\h^0(X,\cO(H)=5$  and $\h^1(X,\cO(H))=4$.   
Indeed, it will turn out that  the component of 
the Hilbert scheme  corresponding to  these surfaces has dimension $38$, hence 
up to projectivities this is a $38-24=14$ dimensional family of abstract
surfaces.
This fits with the fact that the $21$ points have to satisfy a 
condition of codimension $\le 20=4*5$, which leaves us with a
family of collections of points in  $\P^2$ of dimension $\ge 2*21-20=22$.
Up to automorphism of $\P^2$ this leads  to a family of dimension
 $\ge 22-8=14$,
and hence equality holds.
The great difficulty to find points in $\P^2$ in 
very special positions was one of the sources, 
which led Hartshorne to his conjecture.

\noindent

We construct these surfaces via their  
``Beilinson monad": 
Let $V$ be an $n+1$-dimensional vector space over a field $K$ 
and let $W$ be its dual space.  
The basic idea behind a Beilinson monad is to represent 
a given coherent sheaf on $\P^n=\P(W)$ as a homology of a finite
complex
of vector bundles, which are direct sums of exterior powers of the tautological
rank $n$ subbundle $U=\ker(W\tensor \cO_{\P(W)} \to \cO_{\P(W)}(1))$
on $\P(W)$.
(Thus $U \simeq \Omega^1(1)$ is the twisted sheaf of 1-forms. As
Beilinson, we will use the
notation $\Omega^p(p)$ for the exterior powers of $U$.)

The differentials in the monad are given by homogeneous 
matrices over an exterior algebra $E=\bigwedge V$. 
To construct a Beilinson monad for a given coherent sheaf, 
we typically take the following steps:
Determine the type of the Beilinson monad, that is, 
determine the vector bundles of the complex, and then 
find differentials in the monad. 

Let $X$ be a smooth rational surface in $\P^4=\P(W)$ with degree $12$ 
and sectional genus $13$.  
The type of a Beilinson monad for the (suitably twisted) ideal sheaf
of $X$ can be 
derived from the knowledge of its cohomology groups. 
Such information is partially determined from general 
results such as the Riemann-Roch formula and the Kodaira vanishing theorem. 
It is, however, hard to determine the dimensions of all cohomology 
groups needed to determine the type 
of the Beilinson monad. For this reason,  
we assume that the ideal sheaf 
of $X$ has the so-called ``natural cohomology" in some range of
twists. In particular, we assume  that  in each twist $-1\le n \le 6$ 
at most one of the cohomology groups 
$\HH^i(\PP^4,\sI_X(n)$ for $i=0\ldots 4$ is non-zero.
This is an open condition for surfaces in a given component
of the Hilbert scheme.  
Under this assumption the Beilinson monad for the twisted 
ideal sheaf $\sI_X(4)$ of $X$ has the following form:
\begin{eqnarray}
 4\Omega^3(3) \stackrel{A}{\rightarrow} 2\Omega^2(2) \oplus 2\Omega^1(1) 
 \stackrel{B}{\rightarrow} 3\cO. 
 \label{eq:minimal-monad}
\end{eqnarray}
To detect differentials in (\ref{eq:minimal-monad}), we use the following techniques developed recently: (1) the first technique is an exterior algebra method  
due to Eisenbud, Fl{\o}ystad and Schreyer~\cite{efs} 
and (2), the other one is the method using small finite fields and random trials due to Schreyer~\cite{schreyer}.  

\vspace{1mm}
(1) Eisenbud, Fl{\o}ystad and Schreyer  presented 
an explicit version of  the Bernstein-Gel'fand-Gel'fand correspondence. 
This correspondence is an isomorphism 
between the derived category of bounded complexes of  
finitely generated $S$-graded modules and the derived category 
of certain ``Tate resolutions" of $E$-modules, 
where $S=\mathrm{Sym}_K(W)$. 
As an application, they constructed the Beilinson monad from 
the Tate resolution explicitly. 
This enables us to describe the conditions, that 
the differentials in the Beilinson monad must satisfy 
in an exterior algebra context. 

\vspace{1mm}
(2) Let $\M$ be a parameter space for objects in algebraic geometry 
such as the Hilbert scheme or a moduli space. 
Suppose that $\M$ is a subvariety of a rational variety $\GG$ of
codimension $c$.   
Then the probability for a point $p$ in $\GG(\F_q)$ to lie in
$\M(\F_q)$ is  about $(1:q^c)$. 
This approach will be successful if the codimension $c$ is small and 
 the time required to check $ p \not\in \M(\F_q)$ 
is sufficiently small as compared with $q^c$.  
This technique was applied first by Schreyer~\cite{schreyer} 
to find  four different families of smooth surfaces in $\P^4$ 
with degree $11$ and sectional genus $11$ over $\F_3$ by a random search,   
and he provided a method to establish the existence of lifting these surfaces 
to characteristic $0$. 
This technique has been successfully applied to solve various problems in constructive algebraic 
geometry (see \cite{st},  \cite{tonoli} and \cite{bel}).

Singular or Macaulay2 scripts needed to construct and analyse these surfaces 
are available at 
{\tt http://www.math.uni-sb.de/\~{ }ag-schreyer} and  
{\tt http://www.math.colostate.edu/$\sim$abo/programs.html}.

\section{The exterior algebra method}
Our construction of the rational surfaces 
uses the ``Beilinson monad". 
A Beilinson monad represents a given coherent sheaf 
in terms of direct sums of (suitably twisted) bundles of differentials 
and homomorphisms between these bundles, 
which are given by homogeneous matrices over an exterior algebra $E$. 
Recently, Eisenbud, Fl{\o}ystad and Schreyer~\cite{efs} showed 
that for a given sheaf, one can get the Beilinson monad from its ``Tate resolution",    
that is a free resolution over $E$, by a simple functor. 
This enables us to discuss the Beilinson monad in an exterior algebra  context. 
In this section, we take a quick look at the exterior algebra method  
developed by Eisenbud, Fl{\o}ystad and Schreyer. 
\subsection{Tate resolution of a sheaf}
Let $W$ be a $(n+1)$-dimensional vector space over a field $K$, 
let $V$ be its dual space, 
and let $\{x_i\}_{0 \leq i \leq n}$  and 
$\{e_i\}_{0 \leq i \leq n}$ be dual bases of $V$ and $W$ respectively. 
We denote by $S$ the symmetric algebra of $W$  and by $E$ 
the exterior algebra $\bigwedge V$ on $V$. 
Grading on $S$ and $E$ are introduced by $\deg(x)=1$ 
for $x \in W$ and $\deg(e)=-1$ for $e \in V$ respectively. 
The projective space of $1$-quotients of $W$ will be denoted by 
$\P^n=\P(W)$. 

Let $M=\bigoplus_{i\in \integer} M_i$ be a finitely generated $S$-graded module. 
We set 
\[
 \omega_E:=\HHom_K(E,K)=\bigwedge W= E \otimes_K \bigwedge^{n+1} W \simeq E(-n-1)
\]
and 
\[
 F^i:=\HHom_K(E,M_i)\simeq M_i \otimes_K \omega_E. 
\]
The morphism $\phi_i : F^i \rightarrow F^{i+1}$ takes 
the  map $\alpha \in F^i$ to the map
\[
\left[e \mapsto \sum_i x_i\alpha(e_i \wedge e) \right] \in F^{i+1}.  
\] 
Then the sequence 
\[
 \mathbf{R}(M) : \hspace{1mm} \cdots \rightarrow F^{i-1} \stackrel{\phi_{i-1}}{\lra} 
 F^i  \stackrel{\phi_i}{\longrightarrow} F^{i+1} \rightarrow \cdots 
\]
is a complex. This complex is eventually exact. Indeed,  
$\mathbf{R}(M)$ is exact at $\HHom_K(E,M_i)$ for all $i \geq s$ 
if and only if $s>r$,
where $r$ is the Castelnouvo-Mumford regularity of $M$ (see~\cite{efs} for a detailed proof).  
So starting from $\mathbf{T}(M)^{>r}:=\mathbf{T}(M_{> r})$, 
we can construct a doubly infinite exact $E$-free complex $\mathbf{T}(M)$ 
by adjoining a mini\-mal free resolution of the kernel of $\phi_{r+1}$: 
\[
\mathbf{T}(M): \cdots 
\rightarrow T^r \rightarrow T^{r+1}:=\HHom_K(E,M_{r+1}) 
\stackrel{\phi_{r+1}}{\longrightarrow} \HHom_K(E,M_{r+2}) \rightarrow \cdots.  
\]
This $E$-free complex is called the {\it Tate resolution} of $M$. 
Since $\mathbf{T}(M)$ can be constructed by starting from $\mathbf{R}(M_{>s})$, 
$s\geq r$, the Tate resolution depends only on the sheaf 
$\sF = \widetilde{M}$ on $\P(W)$ associated to $M$.
 We call $\mathbf{T}(\sF):=\mathbf{T}(M)$ the {\it Tate resolution} of $\sF$. 
The following theorem gives a description of all the terms of a Tate resolution: 
\begin{theorem}[\cite{efs}]
Let $M$ be a finitely generated graded $S$-module and 
let $\sF := \widetilde{M}$ be the associated sheaf on $\P(W)$.  
Then the term of the complex $\mathbf{T}(\sF)$ with cohomological degree $i$ 
is $\bigoplus_j \HH^j \sF(i-j) \otimes \omega_E$.
\label{th:tate-resolution}
\end{theorem}

Important to us is also the fact the dual complex
$\HHom_E(\mathbf{T}(\sF),E)$ 
stays exact.
 
\subsection{Beilinson monad}
Eisenbud, Fl{\o}ystad and Schreyer~\cite{efs} showed, that 
applying a simple functor to the Tate resolution $\mathbf{T}(\sF)$, 
gives  a finite complex of sheaves whose homology is 
the sheaf $\sF$ itself:  
Given $\mathbf{T}(\sF)$, we define $\Omega(\sF)$ to be the complex 
of vector bundles on $\P(W)$ obtained by replacing each 
summand $\omega_E(i)$ by the  bundle $\Omega^i(i)$. 
The differentials of the complex are given by using isomorphisms 
\[
 \HHom_E(\omega_E(i),\omega_E(j)) \simeq \bigwedge^{i-j} V 
 \simeq \HHom(\Omega^{i}(i),\Omega^j(j)).
\]
\begin{theorem}[\cite{efs}]
 Let $\sF$ be a coherent sheaf on $\P(W)$. 
 Then $\sF$ is the homology of $\Omega(\sF)$ in cohomological 
 degree $0$, and $\Omega(\sF)$ has no homology otherwise. 
\label{th:beilinson-monad}
\end{theorem}

\noindent 
We call $\Omega(\sF)$ the {\it Beilinson monad} for $\sF$. 
\section{Construction }
In this section we will construct 
our family of  rational surfaces $X$ in $\P^4$ with degree $d=12$,  
sectional genus $\pi=13$.   
The construction takes the following  
four steps:
\begin{itemize}
\item[(1)] Analyse  the monad and parts of the Tate resolution.
\item[(2)] 
Find a smooth surface $X$ with the prescribed invariants 
over a finite field of a small characteristic. 
\item[(3)] 
Determine the type of the linear system, 
which embeds $X$ into $\P^4$  to justify that  
the surface $X$ found in the previous step is rational.  
\item[(4)] 
Establish the existence of a lift to characteristic zero. 
\end{itemize}

\subsection{Analysis of the monad and Tate resolution}
Let $K$ be a field, let $W$ be a five-dimensional vector space over $K$ 
with basis $\{x_i\}_{0 \leq i \leq 4}$, 
and let $V$ be its dual space with dual basis $\{e_i\}_{0 \leq i \leq 4}$. 
Let $X$ be a smooth surface in $\P^4=\P(W)$ with the invariants given above. 
The first step is to determine the type of the Beilinson monad for the twisted 
ideal sheaf of $X$, 
which is derived from the partial knowledge of its cohomology groups. 
Such information can be determined from 
general results such as the Riemann-Roch formula and Kodaira vanishing theorem 
(see~\cite{des} for more detail).  
We assume that $X$ has the natural cohomology in the range $-1 \leq j \leq 6$   
of twists:
{
%\fontsize{10pt}{11pt} 
%\selectfont
$$
{
\setlength{\br}{10mm} 
\setlength{\ho}{6mm} 
\fontsize{10pt}{11pt} 
\selectfont 
\begin{xy} 
%
% x-Achse
%
<-3.5\br,0\ho>;<6.5\br,0\ho>**@{-}?>*@{>}?(0.9)*!/^3mm/{j} 
%
% y-Achse 
%
,<-2\br,0\ho>;<-2\br,6\ho>**@{-}?>*@{>}?(0.9)*!/^3mm/{i} 
%
% waagrechte Linien
% 
,0+<-3.5\br,5\ho>;<5.5\br,5\ho>**@{-}
,0+<-3\br,4\ho>;<5\br,4\ho>**@{-}
,0+<-3\br,3\ho>;<5\br,3\ho>**@{-}
,0+<-3\br,2\ho>;<5\br,2\ho>**@{-}
,0+<-3\br,\ho>;<5\br,\ho>**@{-}
%
% senkrechte Linien 
%
,0+<-3\br,0\ho>;<-3\br,5\ho>**@{-} 
,0+<-2\br,0\ho>;<-2\br,5\ho>**@{-} 
,0+<-1\br,0\ho>;<-1\br,5\ho>**@{-} 
,0+<3\br,0\ho>;<3\br,5\ho>**@{-} 
,0+<2\br,0\ho>;<2\br,5\ho>**@{-} 
,0+<1\br,0\ho>;<1\br,5\ho>**@{-}   
,0+<0\br,0\ho>;<0\br,5\ho>**@{-} 
,0+<4\br,0\ho>;<4\br,5\ho>**@{-}
,0+<5\br,0\ho>;<5\br,5\ho>**@{-}  
%
% Eintraege 
%
,0+<-2.5\br,3.5\ho>*{13}
,0+<-.5\br,2.5\ho>*{4}
,0+<.5\br,2.5\ho>*{2}
,0+<1.5\br,1.5\ho>*{2}
,0+<2.5\br,1.5\ho>*{3}
,0+<3.5\br,0.5\ho>*{5}
,0+<4.5\br,0.5\ho>*{29}
,0+<6.5\br,2.5\ho>*{\h^i \sI_X(j)}
%,0+<2.5\br,1.5\ho>*{a}
%,0+<2.5\br,0.5\ho>*{a}
%
% 1.Zeile
%
,<-2.5\br,-.5\ho>*{-1}
,<-1.5\br,-.5\ho>*{0}
,<-.5\br,-.5\ho>*{1}
,<.5\br,-.5\ho>*{2}
,<1.5\br,-.5\ho>*{3}
,<2.5\br,-.5\ho>*{4}
,<3.5\br,-.5\ho>*{5}
,<4.5\br,-.5\ho>*{6}
%
% 1.Spalte
%
,<5.2\br,.5\ho>*{0}
,<5.2\br,1.5\ho>*{1}
,<5.2\br,2.5\ho>*{2}
,<5.2\br,3.5\ho>*{3}
,<5.2\br,4.5\ho>*{4}
\end{xy}
}
$$
}

\noindent
Here a zero is represented by the empty box. 
By Theorem~\ref{th:tate-resolution},
the Tate resolution $\mathbf{T}(\sI_X)[4]=\mathbf{T}(\sI_X(4))$ 
includes an exact $E$-free complex of the following type:
\begin{eqnarray} 
\qquad \rightarrow 4\omega_E(3) 
\rightarrow 2\omega_E(2)\oplus 2\omega_E(1) 
\rightarrow 3\omega_E \oplus 5 \omega_E(-1) \rightarrow 29 \omega_E(-2)
\rightarrow \cdots. 
\label{eq:betti-diagram}
\end{eqnarray}
\noindent
From Theorem~\ref{th:beilinson-monad}, 
it follows therefore, that   the corresponding Beilinson monad for 
$\sI_X(4)$ is of the following type: 
\begin{eqnarray}
0 \rightarrow 4\Omega^3(3) \stackrel{A}{\rightarrow}
 2\Omega^2(2)\oplus2\Omega^1(1) 
\stackrel{B}{\rightarrow} 3\cO \rightarrow 0. 
\label{eq:monad}
\end{eqnarray}

\vspace{2mm}
The next step is to describe what maps $A$ and $B$ 
could be the differentials of the  monad~(\ref{eq:monad}). 
The identifications  
\[
\HHom(\Omega^i(i),\Omega^j(j)) \simeq \HHom_E(\omega_E(i), \omega_E(j)) \simeq \HHom_E(E(i),E(j)), 
\]
allow us to think of the maps $A$ and $B$ as homomorphisms 
between $E$-free modules.  
Since the Tate resolution and its $E$-dual are exact,
the matrix $A$ determines $B$ up to isormorphism. 

However, we start with $B$ in our construction.
To ease our calculations, 
we take the map
$$2\omega_E(1) \stackrel{B_1}{\rightarrow} 3\omega_E $$ 
to be defined by the matrix 
\[
B_1=
\left(
\begin{array}{cc}
e_0 & e_1 \\
e_1 & e_2 \\
e_3 & e_4 
\end{array}
\right), 
\] 
Since the $\GL(5,K)\times \GL(2,K)\times \GL(3,K)$
orbit of this matrix  is dense in $\HHom_E(2\omega_E(1), 3\omega_E)$ 
this is a reasonable mild additional assumption.
The crucial step in the construction is the choice of the map
$$ 3\omega_E \stackrel{C}{\rightarrow}  4\omega_E(-2),$$
where the target $4\omega_E(-2)$ is a free summand of the cokernel
$\Coker(5\omega_E(-1)\to 29\omega_E(-2))$.
Note that $C\circ B = 0$ must hold in the Tate resolution.
The condition $C\circ B_1=0$ means, that 
  $C$ corresponds to a 4-dimensional quotient space
of  
$$T=\Coker(2\Lambda^3 W \stackrel{B_1}{\rightarrow} 3 \Lambda^2 W).$$

An exterior algebra computation proves 
that $\dim T = 10=3*10-2*10$ as expected. Indeed the map to $T$ is
given by the following $10\times 3$ matrix of two forms in $E$:

$$\varphi=
\begin{pmatrix}0&
      0&
      {{e}}_{{3}} {{e}}_{{4}}\\
      0&
      {-{{e}}_{{3}} {{e}}_{{4}}}&
      {{e}}_{{2}} {{e}}_{{3}}-{{e}}_{1} {{e}}_{{4}}\\
      {-{{e}}_{{3}} {{e}}_{{4}}}&
      0&
      {{e}}_{1} {{e}}_{{3}}-{{e}}_{0} {{e}}_{{4}}\\
      0&
      {{e}}_{1} {{e}}_{{4}}-{{e}}_{{2}} {{e}}_{{3}}&
      {{e}}_{1} {{e}}_{{2}}\\
      {{e}}_{{2}} {{e}}_{{3}}-{{e}}_{1} {{e}}_{{4}}&
      {{e}}_{1} {{e}}_{{3}}-{{e}}_{0} {{e}}_{{4}}&
      {-{{e}}_{0} {{e}}_{{2}}}\\
      {{e}}_{0} {{e}}_{{4}}-{{e}}_{1} {{e}}_{{3}}&
      0&
      {{e}}_{0} {{e}}_{1}\\
      0&
      {{e}}_{1} {{e}}_{{2}}&
      0\\
      {{e}}_{1} {{e}}_{{2}}&
      {{e}}_{0} {{e}}_{{2}}&
      0\\
      {{e}}_{0} {{e}}_{{2}}&
      {{e}}_{0} {{e}}_{1}&
      0\\
      {{e}}_{0} {{e}}_{1}&
      0&
      0\\
      \end{pmatrix}$$
Thus we obtain $C$ from a point $[c]\in \GG=\GG(10,4)$ in the
Grassmanian as the product
$C=\varphi\circ c$, where $c\in K^{4\times 10} $ denotes a 
representing $4\times 10$ matrix. 
For these $C$ the condition $C\circ B_1=0$ will be satisfied. 

Consider  
$$\overline \M =\{ [c] \in \GG \mid \exists B_2 \in \HHom(2\omega_E, 3\omega_E(2)) \hbox{ with } C\circ B_2 =0 \}.$$
More precisely,  we consider those $[c]\in \GG$ such that
$$ 0 \to 2\Lambda^4 W \stackrel{B_1}{\rightarrow} 3\Lambda^3 W 
\stackrel{C}{\rightarrow} 4 W \to 0$$
has two dimensional homology in the middle.
The alternating dimensions of the vector spaces
in the complex
add to zero $ 2*5-3*10+4*5=0$. The complex is exact for a 
general choice of $[c]\in \GG$ as we see by a computation in an example.
Thus $[c] \in \GG$, which give the desired  two-dimensional homology 
in the middle,
also give two-dimensional homology at the right. We conclude that 
 $\overline \M\subset \GG$ has codimension at most $4=2*2$ at such points $[c]$.

Once we have  choosen a $[c]\in \overline \M$, we can expect, that $B=(B_1,B_2)$ and $C$ determine the
monad and hence the desired surface, due to the following  Hilbert function 
argument:

The alternating sum of the dimensions in the complex
$$ 0 \to 2\Lambda^3 W\oplus 2\Lambda^4 W \stackrel{B}{\rightarrow} 3\Lambda^2 W 
\stackrel{C}{\rightarrow} 4 K \to 0$$
is $2*10+2*5-3*10+4=4$. Hence we expect a 4 dimensional homology on the right,
which gives the matrix $A$.  

In summary we proved the following proposition.

\begin{proposition}\label{M} There exists a quasi-projective 
subvariety $\M \subset
  \GG(10,4)$ of codimension at most $4$, whose points define a monad of a smooth
  rational surface in $\PP^4$.
The $\PGL(5,\overline K)$ orbit of each  family corresponding to a component
  of $\M$  is an open part of a component of the
  Hilbert scheme of surfaces. 
\end{proposition}

Here $\overline K$ denotes the algebraic
  closure of our ground field $K$.

\begin{proof} Indeed, apart from the condition $[c] \in \overline \M$, all other
  conditions are open conditions.
\end{proof}

However, this does not prove, that $\M$ is non-empty. 
Note  that
$\overline M$ is defined over the integers $\ZZ$.

\subsection{Finite Field search} 

If $\M$ is not empty  we can  expect to find a point in 
$\M(\F_q) \subset \GG(\F_q)$ at a rate of $(1:q^4)$ by Proposition \ref{M}. 
%The procedure
% {\tt numberOfGoodMorForDeg12Surf} 
%counts the number of good points in $\GG(\F_5)$. 
The statistics suggests that there are two  different components of $\overline
\M(\F_5) \subset 
\GG(\F_5)$, whose elements have syzygies with Betti table  
\[
 \begin{tabular}{c|ccccc}
 \hline
 2& 4 &. &. &. \\
 1 &. & 3 & 2 & .\\
 0 &. & . & 2 & 4   \\
 -1 &. & . & . & 5 
 \end{tabular} 
 \quad \mbox{and} \quad 
  \begin{tabular}{c|ccccc}
 \hline
 2& 4 &. &. &. \\
 1 &. & 3 & 2 & .\\
 0 &. & . & 2 & 4   \\
 -1 &. & . & . & 10 
 \end{tabular}  
\]
However, we never obtaiend a Beilinson monad of a surface from an
example with the Betti table of the second type. 
So these points do not belong to $\M(\F_5)$.  
Examples with the first Betti table appeared 
$18$ times in a test of $5^4\cdot 10$ examples. 
It will turn out, that this family has indeed codimension $4$.

\begin{proposition}\label{finite field}
There is a smooth surface in\, $\P^4$ over $\F_5$ with
$d=12$ and $\pi=13$. 
\label{th:small-field}
\end{proposition}

\begin{proof}
By random search, 
%using {\tt randomMonadForDeg12Surf} in Appendix, 
we can find $C\in \M(\F_5)$ and hence $B$ and $A$ satisfying the desired conditions. 
Determine the corresponding maps 
$A: 4\Omega^3(3) \rightarrow 2\Omega^2(2) \oplus 2\Omega^1(1)$ 
and $B=(B_2,B_1) : 2\Omega^2(2) \oplus 2\Omega^1(1) \rightarrow 3\cO$. 
Then compute the homology $\ker(B)/\im(A)$. 
If the homology is isomorphic to the ideal sheaf of a surface 
with the desired invariants, then check smoothness 
of the surface with the Jacobian criterion. If we are lucky, the surface is 
smooth.
If not, we search for a further $C\in \M(\F_5)$.
 
%with {\tt idealOfDeg12Surf}. 
For example the point $[c]\in \overline \M(\F_5)$ represented by the matrix
$$
 c=\begin{pmatrix}{2}&
      {2}&
      {-2}&
      0&
      {-2}&
      {2}&
      {-1}&
      1&
      {-1}&
      {-2}\\
      1&
      {-1}&
      {2}&
      {2}&
      {-1}&
      {2}&
      {2}&
      0&
      {2}&
      {-2}\\
      1&
      {-2}&
      1&
      {-2}&
      0&
      {-1}&
      {-2}&
      {2}&
      1&
      {-2}\\
      {-2}&
      {-1}&
      {-2}&
      {-1}&
      0&
      {2}&
      0&
      {-1}&
      {2}&
      1\\
      \end{pmatrix}$$
leads to a smooth surface in $\P^4$ defined over $\F_5$ of degree $d=12$ and
sectional genus $\pi=13$.
\label{th:first-construction}
\end{proof}
\subsection{Adjunction process}
\label{sec:adjunction_process}
In this subsection, we spot the surface found in the previous step   
within the Enriques-Kodaira classification 
and determine the type of the linear system that embeds $X$ into $\P^4$. 
First of all, 
we recall a result of Sommese and Van de Ven for a surface over $\C$: 
\begin{theorem}[\cite{sv}]
Let $X$ be a smooth surface in $\P^n$ over $\C$ 
with degree $d$, sectional genus  $\pi$, geometric genus $p_g$ and irregularity $q$, let $H$ be its hyperplane class, 
let $K$ be its canonical divisor and let $N=\pi - 1 + p_g -1$. 
Then the adjoint linear system $| H+K |$ defines a birational morphism 
\[
 \Phi=\Phi_{|H+K|} : X \rightarrow \P^{N-1}
\]
onto a smooth surface $X_1$, which blows down precisely all $(-1)$-curves on $X$, 
unless 
\begin{itemize}
 \item[$\mathrm{(i)}$] $X$ is a plane, or Veronese surface of degree $4$, or $X$ is ruled by lines; 
 \item[$\mathrm{(ii)}$] $X$ is a Del Pezzo surface or a conic bundle;
 \item[$\mathrm{(iii)}$] $X$ belongs to one of the following four families: 
 \begin{itemize} 
  \item[$\mathrm{(a)}$]  $X=\P^2(p_1, \dots , p_7)$ embedded by $H \equiv 6L-\sum_{i=0}^7 2E_i$;
  \item[$\mathrm{(b)}$]  $X=\P^2(p_1, \dots , p_8)$ embedded by $H \equiv 6L-\sum_{i=0}^7 2E_i-E_8$;
  \item[$\mathrm{(c)}$]  $X=\P^2(p_1, \dots , p_8)$ embedded by $H \equiv 9L-\sum_{i=0}^8 3E_i$;
  \item[$\mathrm{(d)}$] $X=\P(\sE)$, where $\sE$ is an indecomposable rank $2$ bundle over an
  elliptic curve and $H\equiv B$, where $B$ is a section $B^2=1$ on $X$. 
 \end{itemize}
\end{itemize}
\label{th:adjoint-linear-system}
\end{theorem}
\begin{proof}
 See~\cite{sv} for the proof. 
\end{proof}

\noindent 
Setting $X=X_1$ and performing the same operation repeatedly, 
we obtain a sequence 
\[
 X \rightarrow X_1 \rightarrow X_2 \rightarrow \cdots \rightarrow X_k. 
\] 
This process will be terminated if $N-1 \leq 0$. 
For a surface with nonnegative Kodaira dimension, 
one obtains the minimal model at the end of the adjunction process. 
If the Kodaira dimension equals $-\infty$, 
we end up with a ruled surface, 
a conic bundle, a Del Pezzo surface, 
$\P^2$, or one of the few exceptions of Sommese and Van de Ven. 

It is not known, whether the adjunction theory holds over a finite field. 
However, we have the following proposition:
\begin{proposition}[\cite{ds}, Prop. 8.3]
Let $X$ be a smooth surface over a field of arbitrary  characteristic. 
Suppose that the adjoint linear system $|H+K|$ is base point free. 
If the image $X_1$ in $\P^N$ under the adjunction map $\Phi_{|H+K|}$ 
is a surface of the expected degree $(H+K)^2$, the expected sectional genus 
$\frac{1}{2}(H+K)(H+2K)+1$ and  with $\chi(\cO_X)=\chi(\cO_{X_1})$, 
then $X_1$ is smooth and $\Phi: X\rightarrow X_1$ is a simultaneous 
blow down of the $K_1^2-K^2$ many exceptional lines on $X$. \qed
\label{th:adjunction-p}
\end{proposition}

\begin{remark}
The union of the exceptional divisors contracted in each step is defined over 
the base field.
\end{remark}
In~\cite{des} and~\cite{ds}, it is described 
how to compute the adjunction process for a smooth surface given by explicit equations 
(see~\cite{ds} for the computational details).   
Let $X$ be the smooth surface found in the previous step. 
The computation for the adjunction process in characteristic $5$ gives 
\begin{eqnarray}
H\equiv 12L-\sum_{i_1=1}^2 4E_{i_1}-\sum_{i_2=3}^{11} 3E_{i_2} 
-\sum_{i_3=12}^{14} 2E_{i_3}-\sum_{i_4=15}^{21} E_{i_4}, 
\label{eq:ample-divisor}
\end{eqnarray}
where $L$ is the class of a line in $\P^2$. 
This process ends with a Del Pezzo surface of degree $7$, 
which is the blowing up of $\P^2$ in two points. 
Therefore we can conclude that $X$ is rational. 

%\begin{remark}
%Another way to prove rationality of $X$ is to count the number of 
%$6$-secant lines to $X$.
%First we prove the following claim: 
%Let $Y$ be a smooth surface in $\P^4$ with $d=12$, $\pi=13$ and 
%lines to $Y$. Then $Y$ is rational. 

%To prove this claim, we use 
%the Le Barz's formula $N_6(d,\pi, \chi)$. This formula gives us $N_6(12,13,1)=8$,  
%that equals the number 
%of $6$-secant lines to $Y$ plus the number of exceptional lines 
%on $Y$, if there are at most a finite number of $6$-secant lines to $Y$ and 
%if there are no lines on $Y$ with nonnegative self-intersection number~\cite{leBarz}. 

%Let $H_2$ be the hyperplane class of the second adjoint surface 
%$Y_2$ and let $K_2$ be the canonical divisor of $Y_2$. 
%Then we have $H_2 \cdot K_2= -12+a$, 
%where $a$ is the number of exceptional lines on $Y$. 
%The Le Barz's formula tells us that $Y$ can have at most eight exceptional lines 
%if there are a finite number of $6$-secant lines to $Y$. 
%So $H_2 \cdot K_2 < 0$, and hence $Y$ is rational.  

%Next we show that $X$ has only one $6$-secant line. 
%The union of $6$-secant lines to $X$ is contained in 
%all the quintics that contain $X$. 
%With Computer algebra  we can check that 
%\[
% V(\HH^0 \sI_X(5))=X \cup L_0, 
%\]
%where $L_0$ is a line. 
%So $L_0$ is the only $6$-secant line to $X$, 
%and hence $X$ is rational. 
%\end{remark}

%
\subsection{Lift to characteristic $0$} 
\label{lift}
In the previous step, we constructed a smooth surface in $\P^4$
defined over $\F_5$. 
However, our main interest is the field of complex numbers $\C$.  
In this section,  
we show the existence of a lift to characteristic $0$ as follows: 
Let $\M$ and $\GG$ be given as in the previous subsections.  
\begin{proposition}[\cite{schreyer}]
Let $[c] \in \M(\F_p)$ be a point, 
where $\M\subset \GG$ has codimension $4$. 
Then there exist  a number field $\LL$, 
a prime $\mathfrak{p}$ in $\LL$  with residue field 
$\cO_{\LL,\mathfrak{p}} /\mathfrak{p}\cO_{\LL,\mathfrak{p}} \simeq
\F_p$ and a family of surfaces $\sX$ defined over
$\cO_{\LL,\mathfrak{p}}$
with special fiber the surface $X$ defined over $\F_p$ corresponding to $[c]$.
Furthermore, since the surface $X/\F_p$ corresponding to $[c]$ is smooth, 
the surface $X/\LL$ corresponding to the generic point 
of $\Spec \LL \subset \Spec \cO_{\LL,\mathfrak{p}}$ is also smooth. 
\label{th:lift}
\end{proposition}
\begin{proof}
Let $p$ be a prime number. 
If this is not the case, $\Z$ has to be replaced by the ring of integers 
in a number field which has $\F_p$ as the residue field. 

Since $\M$ has pure codimension $4$ in $[c]$, 
there are four hyperplanes $H_1,\dots,H_4$ in $\GG$, 
such that $[c]$ is an isolated point of $\M(\overline{\F}_p) \cap H_1 \cap \cdots \cap H_4$. 
We may assume that $H_1, \dots, H_4$ are defined over $\Spec\Z$ 
and that they meet transversally in $[c]$. 
This allows us to think that $\M \cap H_1 \cap \cdots \cap H_4$ 
is defined over $\Z$. Let $Z$ be an irreducible component of $\M_{\Z}$ 
containing $C$.
Then $\dim Z=1$. 

The residue class field of the generic point of $Z$ 
is a number field $\LL$ that is 
finitely generated over $\Q$, 
because $\M$ is projective over $\Z$. 
Let $\cO_{\LL}$ be the ring of integers of $\LL$ 
and let $\mathfrak{p}$ be a prime ideal corresponding to $[c] \in Z$. 
Then $\Spec\cO_{\LL, \mathfrak{p}} \rightarrow Z \subset \M$ 
is an $\cO_{\LL,\mathfrak{p}}$-valued point which lifts $[c]$.

Performing the construction of the surface over $\cO_{\LL,\mathfrak{p}}$ 
gives a flat family $\sX$ of surfaces over $\cO_{\LL,\mathfrak{p}}$. 
Since smoothness is an open property, and since the special fiber 
$X=\sX_{\mathfrak{p}}$ is smooth, 
the general fiber $\sX_{\LL}$ is also smooth. 
\end{proof}
Next, we argue that the adjunction process of the surface over the number field 
$\LL$ has the same numerical behavior: 
\begin{proposition}[\cite{ds}, Cor. 8.4]
Let $\sX \rightarrow \Spec \cO_{\LL,\mathfrak{p}}$ be a family as in 
Proposition~\ref{th:lift}. 
If the Hilbert polynomial of the first adjoint surface of $X=\sX \tensor \F_q$ 
is as expected, and if 
$\HH^1(X,\cO_X(-1))=0$, then 
the adjunction map of the general fiber $\sX_{\LL}$ blows down 
the same number of exceptional lines as the adjunction map of the special fiber $X$. 
\label{th:adjunction-0}\qed
\end{proposition}

{\it Last step in the proof of Theorem \ref{main}.} 
Let $[c]$ be the element of $\M(\F_5)$, which gives
the surface in Proposition \ref{finite field}.   
We check,  that $[c]$ satisfies the condition of Proposition \ref{th:lift} 
by computing the Zariski tangent space $T_{\M,[c]}$ at $[c]$. Our
computation shows that $\codim T_{\M,[c]} =4$. So $\M$ is smooth of codimension
$4$ at $[c]$, and $[c]$ and hence the surface lift to a number field.

Finally we count dimension. Our component $\M \subset \GG(10,4)$ containing
$[c]$ has
codimension $4$, hence dimension $4*(10-4)-4=20$. The
normalization of $B_1$ (up to conjugation) gives additional $18$ parameters,
because the Hilbert scheme of cubic scrolls in $\P(V)$ has dimension $18$.
So the component of the Hilbert scheme, that contains our surface, has 
dimension $38$. 
\qed

%
%%%%%%%%%%%%%%%%%%%%%%%%%%%%%%%%%%%%%%%%%%%%%%%%%%%%%%%%%%%%%%%%%%%%%%%%%%%%
%      At the end of the article: acknowledgements and literature 
%
\section*{Acknowledgements}
%

%
%     INSERT here acknowledgements 
%

%


\begin{thebibliography}{99}
%
%   References should be in the format [1], [2], ...
%   
% \bibitem{xxx} Name, A.: {\it Title of Article.} Journal {\bf 1}, 333--345
% (2005). 
%
%\bibitem{abo}H.~Abo,  
%Macaulay 2 scripts for finding rational surfaces in $\P^4$ with degree 12.  
%Available at {\tt http://www.math.colostate.edu/$\sim$abo/programs.html}. 
%
%
%\bibitem{ar}H.~Abo and K.~Ranestad,  
%{\it Construction of rational surfaces in projective fourspace},   
%in preparation.  
%
\bibitem{bel}H.~C. Bothmer, C.~Erdenberger and K.~Ludwig. 
{\it A new family of rational surfaces in $\P^4$},
math.AG/0404492, 2004 
%

%
\bibitem{des}W.~Decker, L.~Ein and F.-O.~Schreyer, 
{\it Construction of surfaces in $\P^4$},
J.~Algebr.~Geom. 
{\bf 2}
(1993)
185--237
%
%
\bibitem{de}W.~Decker and D.~Eisenbud,  
{\it Sheaf algorithms using the exterior algebra},
In: D.~Eisenbud, D.R.~Grayson, M.E.~Stillman, B.~Sturmfels (Eds), 
Computations in Algebraic Geometry with {\it Macaylay 2},   
Algorithms~Comput.~Math. 
{\bf 8} 
Springer, 
Berlin
(2002)
215--249
%
%
\bibitem{ds}W.~Decker and F.-O.~Schreyer,  
{\it Non-general type surface in $\P^4$: some remarks on bounds and constructions}, 
J.~Symbolic Comput., 
{\bf 25}  
(2000)  
545--582 
%
%
\bibitem{efs} D.~Eisenbud, G.~Fl{\o}ystad and F.-O.~Schreyer,   
{\it Sheaf cohomology and free resolutions over exterior algebras}, 
Trans. Amer. Math. Soc. 
{\bf 355}  
(2003)  
4397--4426 
%
\bibitem{ep} G.~Ellingsrud and C.~Peskine 
{\it Sur les surfaces lisse de $\P^4$},
Invent. Math. 
{\bf 95}
(1989)
1--12
%
\bibitem{gs}D.~Grayson and M.~Stillman, 
(1991) 
Macaulay 2, a software system for research in algebraic geometry. 
Available at {\tt http://www.math.uiuc.edu/Macaulay2}. 
%
\bibitem{Sing} G.-M.~Greuel, G.~Pfister and H.~Sch\"onemann,
{\sc Singular} 2.0,
A Computer Algebra System for Polynomial Computations,
Centre for Computer Algebra, University of Kaiserslautern,
{\tt http://www.singular.uni-kl.de}
%
%\bibitem{leBarz}P.~Le Barz, 
%{\it Formules pour les multisecantes des surfaces},
%C.~R.~Acad.~Sc.~Paris,   
%{\bf 292}
%(1981)
%797--799
%
\bibitem{schreyer}F.-O.~Schreyer,    
{\it Small fields in constructive algebraic geometry}, 
Lecture Notes in Pure and Appl. Math.,  
{\bf 179}
(1996)
221--228.
%
\bibitem{st}F.-O.~Schreyer and F.~Tonoli,    
{\it Needles in a Haystack: Special varieties via small field}, 
In: D.~Eisenbud, D.R.~Grayson, M.E.~Stillman, B.~Sturmfels (Eds), 
Computations in Algebraic Geometry with {\it Macauylay 2},   
Algorithms~Comput.~Math. 
{\bf 8} 
Springer, 
Berlin
(2002)
215--249
%
\bibitem{sv} A.~J. Sommese and A.~Van de Ven,  
{\it On the adjunction mapping}, 
Math. Ann. 
{\bf 278}
(1987)
593--603.
%

\bibitem{tonoli}F.~Tonoli,  
{\it Construction of Calabi-Yau $3$-folds in $\P^6$},
J. Algebr. Geom 
{\bf 13}
(2004)
209--232.
 
\end{thebibliography}
\end{document}